\documentclass[11pt]{amsart}
\usepackage{amsmath, amscd, amsthm, amsfonts, amssymb, mathrsfs}
\usepackage{epic,eepic,epsf,epsfig}
\usepackage{amsfonts,srcltx,mathrsfs}
\usepackage{amsmath}
\usepackage{amssymb}
\usepackage{amsbsy}
\usepackage{float}
\usepackage{graphicx}
\usepackage{amsfonts}
\usepackage{color}
\usepackage{tikz}
\usepackage{comment}

\usepackage{url}
\newtheorem{lem}{Lemma}[section]
\newtheorem{theorem}[lem]{Theorem}

\newtheorem{case}{Case}
\newtheorem{prop}[lem]{Proposition}
\newtheorem{rem}[lem]{Remark}

 \def\b{\beta}   
\def\r{\rho} \def\s{\sigma}   
 \def\ld{\lambda}

 \def\lg{\langle} \def\rg{\rangle}
\def\nd{\mathrel{\bigm|\kern-.7em/}}

\def\f{\noindent}

\def\PSL{\hbox{\rm PSL}}

\def\U{\hbox{\rm PSU}}
\def\Aut{\hbox{\rm Aut}}

\def\Syl{\hbox{\rm Syl}}

\def\Aut{\hbox{\rm Aut}}

\def\mod{\hbox{\rm mod }}

\def\demo{\f {\bf Proof.}\hskip5pt}

\def\mz{{\mathbb Z}}
\def\P{\mathcal {P}}

\begin{document}

\title{Chiral polytopes of order $2p^m$}

\author{Ting-Ting Kong}
\address{Ting-Ting Kong, Department of Mathematics, Shanxi Key Laboratory of Cryptography and Date Security, Shanxi Normal University, TaiYuan,
030031, P.R. China}
\email{703537@sxnu.edu.cn}

\author{Yan-Quan Feng}
\address{Yan-Quan Feng, Department of Mathematics, Beijing Jiaotong University, Beijing,
100044, P.R. China}
\email{yqfeng@bjtu.edu.cn}

\author{Dong-Dong Hou}
\address{Dong-Dong Hou, Department of Mathematics, Shanxi Normal University, TaiYuan,
030031, P.R. China}
\email{holderhandsome@bjtu.edu.cn}

\author{Dimitri Leemans}
\address{Dimitri Leemans, D\'epartement de Math\'ematique, Universit\'e Libre de Bruxelles, C.P.216, Boulevard du Triomphe, 1050 Bruxelles Belgium}
\email{leemans.dimitri@ulb.be}

\author{Hai-Peng Qu}
\address{Hai-Peng Qu, Department of Mathematics, Shanxi Normal University, TaiYuan,
030031, P.R. China}
\email{orcawhale@163.com}

\date{}

\begin{abstract}
Let $\mathcal{P}$ be a chiral polytope with type $\{k_1, k_2\}$ and $G=\Aut(\mathcal{P})$.
Suppose $|G|=2p^m$, where $k_1, k_2\geq 3$ and $p$ is an odd prime.
Let $P$ be a Sylow $p$-subgroup of $G$. We prove that $G \cong P \rtimes \mathbb{Z}_2$, $d(P)=2$, $P' \ne 1$(so $m \geq 3$)
and up to duality, $\{k_1, k_2\}=\{p^{l_1}, 2p^{l_2}\}$ for some integral $l_1, l_2 \geq 1$. Moreover, 
we show that $\mathcal{P}$ is tight $(k_1k_2=2p^m$) if and only if $P$ is metacyclic group. Furthermore, if $m=3$ or $4$, then
$\mathcal{P}$ must be tight, and if $m \geq 5$, where either $m$ is odd, or $m$ is even and $m \geq p+3$, 
there exists a non-tight chiral polytope $\mathcal{P}$.

\end{abstract}
\maketitle

\section{Introduction}\label{sec1}

Abstract polytopes are combinatorial structures with properties that generalise those of
classical polytopes.  In many ways they are more fascinating than convex polytopes
and tessellations.  Highly symmetric examples of abstract polytopes include not
only classical regular polytopes such as the well known Platonic solids, and more exotic structures such as
the $120$-cell and $600$-cell, but also regular maps on surfaces (such as Klein's quartic);
see~\cite[Chapter 8]{BooksCoxeter} for example.

Roughly speaking, an abstract polytope is a partially ordered set $(\P, \leq)$ (or $\P$ for short) endowed with
a rank function, satisfying certain conditions that arise naturally from a geometric
setting (see Section~\ref{definitions} for definitions). Such objects were proposed by Branko Gr\"unbaum in the seventies; their definition
(initially as `incidence polytopes') and theory were developed by Ludwig Danzer and Egon Schulte. They are closely related to `polyhedral geometries' that Jacques Tits defined in 1961~\cite{Tits1961}. 

An \emph{automorphism} of an abstract polytope $\P$ is an order-preserving permutation
of the elements of $\P$, and every automorphism of $\P$ is uniquely determined by its action
on any maximal chain in $\P$ (known as a `flag' of $\P$). The
most symmetric examples of abstract polytopes are the regular ones, where the automorphism group acts transitively on the set of flags.
The comprehensive book written by Peter McMullen and Egon Schulte~\cite{ARP} is considered the main reference on this subject.

Another interesting class of examples is the one of chiral polytopes.
For these polytopes, their automorphism group has two orbits on the set of flags, with the extra property that any two flags that
differ in a single element lie in different orbits.
The study of chiral abstract polytopes
was pioneered by Schulte and Asia Ivic Weiss (see~\cite{ChiralPolytopes, ChiralPolytopesWithPSL} for example).
For quite some time, the only known finite examples
of chiral polytopes had ranks 3 and 4, but then some finite examples of rank 5 were
constructed by Marston Conder, Isabel Hubard and Tomaz Pisanski~\cite{2008ChiralPolytopesOfSmallRank}. Daniel Pellicer was the first to show that chiral polytopes of arbitrarily large rank exist~\cite{Pellicer2010}. 

It is a natural question to try to describe all pairs $(\P, G)$, where $\P$ is a chiral polytope and $G$ is the automorphism group of $\P$. 
Conder computed all chiral polytopes with automorphism group of order at most 1000 in 2012\footnote{These computations have recently been pushed to cover all polytopes with up to 4000 flags~\cite{atlas1}.}.
Roughly at the same time, Michael Hartley, Hubard and Dimitri Leemans built in~\cite{atlas2} an atlas of all chiral polytopes whose automorphism 
group is an almost simple group $G$ such that $S  \leq G \leq \Aut(S)$ and $S$ is a simple group of order less than 1 million appearing in the Atlas of Finite Groups by Conway et al. Recently, Primoz Potočnik, Pablo Spiga and Gabriel Verret computed all chiral 3-polytopes with automorphism group of order at most 12000~\cite{GroupOrder6000}.

An interesting case consists of pairs $(\mathcal{P}, G)$ with $G$ simple or almost simple.
In 2015, Conder, Hubard, Pellicer and Eugenia O'Reilly-Regueiro~\cite{ChiralPolytopesOfANSN} provided a construction of chiral 4-polytopes for $A_n$ and $S_n$. After ten years, Conder, Hubard and O'Reilly-Regueiro extended their initial idea and proved that for every $d \geq 5$, there exists a chiral $d$-polytope with automorphism group isomorphic to $A_n$, and a chiral $d$-polytope with  automorphism group isomorphic to $S_n$ for all but finitely many $n$, see~\cite{HighRankChiralPolytopesOfANSN}. Also, Wei-Juan Zhang gave some new construction for chiral $4$-polytopes 
with symmetric automorphism groups~\cite{2023+Zhang+ChiralPolytopes}.
Hubard and Leemans dealt with the Suzuki simple groups in~\cite{2014+Hubard+Leemans}. In a series of papers, the case where the automorphism group of a chiral polyhedra is an almost simple group with socle $\PSL(2,q)$ has been studied~\cite{Conder2008,2015+Elliott+Leemans+O'Reilly,2016+Leemans+Moerenhout,2017+Leemans+Moerenhout+O'Reilly}. Moreover, in 2022, Leemans and Adrien Vandenschrick constructed chiral 4-polytopes with automorphism group $\PSL(3, q)$ for every prime power $q \geq 5$~\cite{2022+Leemans+Vandenschrick}. 
For chiral $3$-polytopes, Conder, Potočnik and Jozef \v{S}ir{\'a}\v{n} proved that the groups $\PSL(2,q)$ cannot be automorphism groups of such polyhedra~\cite{Conder2008}. Leemans and Martin Liebeck proved that every finite simple group acts as the automorphism group of a chiral $3$-polytope, apart from
the groups $\PSL(2, q)$, $A_7$, and maybe also $\PSL(3,q), \U(3,q)$ and ~\cite{2017+Leemans+Liebeck}. Antonio Breda and Domenico Catanalo then showed that $\PSL(3,q)$ and $\U(3,q)$  cannot be automorphism groups of chiral polyhedra (with the help of Leemans and Liebeck for the groups $\U(3,q)$)~\cite{BC2021}.

Another interesting case consists of the pairs $(\mathcal{P}, G)$ with $G$ a soluble group.
When it comes to soluble groups, except for abelian groups,
the first families that come to mind are groups of orders $2^n$ or $2^np^m$, where $p$ is an odd prime. Schulte and Weiss mentioned in~\cite{Problem} that the order $2^n$ and $2^np$ prove to be more difficult than others and suggested to try to classify chiral polytopes with these orders~\cite[Problem 30]{Problem}.
Zhang constructed in her PhD thesis~\cite{Zhang'sPhdThesis} a number of chiral polytopes
of types $\{4, 4, 4\}$, $\{4, 4, 4, 4\}$ and $\{4, 4, 4, 4, 4\}$, with automorphism groups of orders $2^{10}, 2^{11}, 2^{12}$,
and $2^{15}, 2^{16}, \cdots, 2^{22}$, and $2^{18}, 2^{19}$, respectively. 
Recently, Zhang presented a general construction for chiral polytopes of type $\{4, 4, \cdots, 4\}$ with rank $4$, $5$ and $6$, which are
obtained as boolean covers of the unique tight regular polytope of the same type in~\cite{2023+Zhang+ChiralPolytopeswithtype444}.
Conder, Yanquan Feng and Dong Dong Hou gave two infinite families of chiral 4-polytopes of type $\{4, 4, 4\}$ with solvable automorphism groups~\cite{2021+Conder+Feng+Hou}.
In 2024, Hou, T. T. Zheng and R. R. Guo proved that for $n \geq 9$, there exists a chiral 3-polytope of type $\{4, 8\}$ with automorphism group of order $2^n$~\cite{2024+Hou+Zheng+Guo}.
Gabe Cunningham gave some infinite families of tight chiral $3$-polytopes of type $\{k_1, k_2\}$ with automorphism group of order $2k_1k_2$~\cite{TightChiralPolyhedra}.

Note that if a chiral polytope $\mathcal{P}$ has rank at least three, the order of its automorphism group $G$ must be even.
In this paper, we focus on groups $G$ acting on chiral polytopes of rank three, with $|G|=2p^m$, $p$ an odd prime, and $m$ a strictly positive integer.
We first prove the following theorem (where $d(G)$ is the size of a smallest generating set of $G$).

\begin{theorem}\label{maintheorem1}
Let $p$ be an odd prime and $m$ be a positive integer. 
Let $G$ be a group of order $2p^m$.
Let $(G,\{\sigma_1,\sigma_2\})$ be the automorphism 
group of a chiral polyhedron $\mathcal{P}$ of type $\{k_1, k_2\}$.
Then
\begin{itemize}
\item [\rm {(1)}] $G = P \rtimes \lg \sigma_1\sigma_2\rg$ with $P \in \Syl_p(G)$;
\item [\rm {(2)}] Up to duality, $\{k_1, k_2\}=\{p^{l_1}, 2p^{l_2}\}$ for some positive integers $l_1, l_2$;
\item [\rm {(3)}] $P=\lg \sigma_1, \sigma_2^2\rg, d(P)=2$ and $P$ is not abelian (so $m \geq 3$);
\item [\rm {(4)}] $\mathcal{P}$ is tight, that is $k_1k_2=2p^m$, if and only if $P$ is a metacyclic group.
\end{itemize}
\end{theorem}

We then prove the following theorem, showing that if $m =3$ or $4$, then $\P$ is a tight polyhedron.

\begin{theorem}\label{maintheorem2}
Let $p$ be an odd prime and $m = 3$ or $4$.
Let $G$ be a group of order $2p^m$.
let $(G,\{\sigma_1,\sigma_2\})$ be the automorphism 
group of a chiral polyhedron $\mathcal{P}$ of type $\{k_1, k_2\}$.
Then $\mathcal{P}$ is tight.
\end{theorem}

\begin{rem}
As a by-product of the proof of Theorem~\ref{maintheorem2}, we obtain the following families of non-tight regular polyhedra:

\begin{itemize}
\item regular polyhedra of type $\{p, 2p\}$ of order $4p^3$;
\item regular polyhedra of type $\{p, 2p\}$ of order $4p^4$;
\item regular polyhedra of type $\{p^2, 2p\}$ of order $4p^4$;
\item regular polyhedra of type $\{p, 2p^2\}$ of order $4p^4$.
\end{itemize}
In particular, the first class and part of the second class of 
these polyhedra are identical to those Feng, Hou, Leemans and Qu constructed in~\cite{2025+Hou+Feng+Leemans+Qu}.
\end{rem}

Cunningham proved in~\cite{TightChiralPolyhedra} that there exists a tight chiral polyhedron $\mathcal{P}$ of type $\{p^{l_1}, 2p^{l_2}\}$ for $l_1 \geq 2, l_2 \geq 1$ and $l_1 \geq l_2$.
In the case where $\mathcal{P}$ is not tight, we prove the following theorem. The proof of that theorem is constructive, meaning that we give an explicit way to construct such a polyhedron.

\begin{theorem}\label{maintheorem3}
Let p be an odd prime. Let $m \geq 5$ be an integer such that either $m$ is odd, or $m$ is even and $m \geq p+3$.  
Then there exists a non-tight chiral poyhedron $\mathcal{P}$ with $|\Aut(\mathcal{P})|=2p^m$.
\end{theorem}

Let us point out that, as in the regular case (see~\cite[Section 1]{2025+Hou+Feng+Leemans+Qu}) chiral polytopes constructed from soluble groups are extremely rare.
For example, there are 49,910,526,325 groups of order $\leq 2000$~\cite{Groups2000}. These groups are readily accessible in {\sc Magma}~\cite{BCP97} except for those of order 1024. Even though those of order 1024 are not in {\sc Magma}, we know they are all soluble.
One can check that out of all groups of order at most 2000, there are 49,910,525,301  that are soluble groups and 1024 that are non-soluble.
According to the data collected by Conder~\cite{atlas1}, the number of chiral polyhedra having an automorphism group of order at most 2000 that is soluble (resp. non-soluble) is 2665 (resp. 165).
Hence the ratio (number of chiral group representations/number of groups) is 0.000000053\% for soluble groups and 16.11\% for non-soluble groups.
This shows that it is really not easy to find soluble groups that are the automorphism group of some chiral polytopes. More importantly, if we refine by rank, we get, in the soluble case, 2389 polytopes of rank 3, 276 of rank 4 and none of rank 5. In the non-soluble case, we get 39 polytopes of rank 3, 123 of rank 4 and three of rank 5.

The paper is organised as follows.
In Section~\ref{sec2}, we give the necessary background to understand this paper.
In Section~\ref{sec3}, we give several lemmas about finite $p$-groups of maximal class that are useful in this paper. 
In Section~\ref{sec4}, we prove Theorem~\ref{maintheorem1}.
In Section~\ref{sec5}, we prove Theorem~\ref{maintheorem2}.
Finally, in Section~\ref{sec6}, we prove Theorem~\ref{maintheorem3}.

\section{Preliminaries}\label{sec2}

In this section we give some further background that may be helpful for the rest of the paper (see~\cite{GroupBookss,Huppert,ARP,ChiralPolytopes,XQ} for more details).

\subsection{Abstract polytopes: definition, structure and properties}\label{definitions}

An abstract polytope of rank $n$ is a partially ordered set $(\mathcal{P}, \leq)$ endowed with a strictly monotone rank function
with range $\{-1, 0, \cdots, n\}$, which satisfies four properties (P1) to (P4) described below.

The elements of $\mathcal{P}$ are called \emph{faces} of $\mathcal{P}$. More specifically, the elements of $\mathcal{P}$ of rank $j$ are called {\em $j$-faces},
and a typical $j$-face is denoted by $F_j$.
Two faces $F$ and $G$ of $\mathcal{P}$ are said to be \emph{incident} with each other if $F \leq G$ or $F \geq G$ in $\mathcal{P}$.
A \emph{chain} of $\mathcal{P}$ is a totally ordered subset of $\mathcal{P}$, and is said to have \emph{length} $i$ if it contains exactly $i+1$ faces.
The maximal chains in $\P$ are called the \emph{flags} of $\mathcal{P}$.
Two flags are said to be $j$-\emph{adjacent} if they differ in just one face of rank $j$,
or simply \emph{adjacent} (to each other) if they are $j$-adjacent for some~$j$.
Also if $F$ and $G$ are faces of $\mathcal{P}$ with $F \leq G$, then the set $\{\,H \in \mathcal{P} \ |\ F \leq H \leq G\,\}$ is
called a \emph{section} of $\mathcal{P}$, and is denoted by $G/F$.  Such a section has rank $m-k-1$,
where $m$ and $k$ are the ranks of $G$ and $F$ respectively.
A section of rank $d$ is called a $d$-section.

\smallskip
We can now give the four conditions that are required of $\P$ to make it an abstract polytope.
\begin{itemize}
\item [(P1)]  $\P$ contains a least face and a greatest face, denoted by $F_{-1}$ and $F_n$, respectively.
\item [(P2)]  Each flag of $\P$ has length $n+1$ (so has exactly $n+2$ faces, including $F_{-1}$ and $F_n$).
\item [(P3)]  $\P$ is \emph{strongly flag-connected}, which means that any two flags $\Phi$ and $\Psi$ of $\P$ can
be joined by a sequence of successively adjacent flags $\Phi=\Phi_{0}, \Phi_1, \cdots, \Phi_k=\Psi$, each of which contains $\Phi \cap \Psi$.
\item [(P4)]  The rank $1$ sections of $\P$ have a certain homogeneity property known as the \emph{diamond condition}: if $F$ and $G$ are incident faces of $\P$, of ranks $i-1$ and $i+1$, respectively, where $0 \le i \le n-1$,
then there exist precisely \emph{two} $i$-faces $H$ in $\P$ such that $F< H< G$.
\end{itemize}

An easy case of the diamond condition occurs for polytopes of rank 3 (or polyhedra). In that case, 0-faces are called {\em vertices}, 1-faces are called {\em edges} and 2-faces are called {\em faces}. If $v$ is a vertex of some face $f$,
then there are two edges that are incident with both $v$ and $f$.

\smallskip
If $F_{n-1}$ is a face of rank $n-1$, then
the section $F_{n-1}/F_{-1}$ is also called a {\em facet} of $\P$,
while if $F_0$ is a face of rank 0, then the section
$F_{n}/F_{0} = \{G \in \P \ |\ F_{0} \leq G\}$ is called a {\em vertex-figure} of $\P$ at $F_0$.
Every $2$-section $G/F$ of $\P$ is isomorphic to the face lattice of a polygon.
If it happens that the number of sides of every such polygon depends only on the rank of $G$,
and not on $F$ or $G$ itself, then we say that the polytope $\P$ is {\em equivelar}.
In this case, if $k_i$ is the number of edges of a $2$-section between an $(i-2)$-face and an $(i+1)$-face of $\P$,
for $1 \leq i \leq n$, then the ordered set $\{k_1, k_2, \cdots, k_{n-1}\}$ is called the {\em Schl\"afli type} of $\P$.
(For example, if $\P$ has rank 3, then $k_1$ and $k_2$ are respectively the size of each face and the valency of each vertex.)

\subsection{Automorphisms of polytopes}

An \emph{automorphism} of an abstract polytope $\P$ is an order-preserving permutation of its elements.
In particular, every automorphism preserves the set of faces of any given rank.
Under permutation composition, the set of all automorphisms of $\P$ forms a group, called the {\em automorphism
group} of $\P$, and denoted by $\Aut(\P)$ or sometimes more simply as $G$.
Also it is not difficult to use the diamond condition and strong flag-connectedness to prove that if an automorphism
preserves one flag of $\P$, then it fixes every flag of $\P$ and hence every element of $\P$.
It follows that $G$ acts semi-regularly on flags of $\P$.

\smallskip
A polytope $\P$ is said to be \emph{regular} if its automorphism group $G$ acts transitively
(and hence regularly) on the set of flags of $\P$.
In this case, the number of automorphisms of $\P$ is as large as possible, and equal to the number of flags of $\P$.
In particular, $\P$ is equivelar, and the stabiliser in $G$ of every 2-section of $\P$ induces the full dihedral group
on the corresponding polygon.
Moreover, for a given flag $\Phi$ and for every $i \in \{0,1,\dots,n-1\}$, the polytope $\P$ has a unique
automorphism $\rho_i$ that takes $\Phi$ to its unique $i$-adjacent flag $\Phi^{i}$.
The automorphisms $\rho_0,\rho_1,\dots,\rho_{n-1}$ generate $G$ and satisfy the defining
relations for the string Coxeter group $[k_1, k_2, \cdots, k_{n-1}]$, where the $k_i$ are as given in the
previous subsection for the Schl\"afli type of $\P$.
Here, the {\em string Coxeter group} $[k_1, k_2, \cdots, k_{n-1}]$, is defined as the group with presentation

$$\lg \r_0, \r_1, \cdots, \r_{n-1} \ |\ \r_i^2=1 \ {\rm for} \ 0 \leq i \leq n-1, (\r_i\r_{i+1})^{k_{i+1}}=1 \ {\rm for} \ 0 \leq i \leq n-2,$$
$$(\r_i\r_j)^2=1 \ {\rm for} \ 0 \leq i < j-1  < n-1\rg.$$

They also satisfy an intersection condition, which follows from the diamond and strong flag-connectedness
conditions.
These and many more properties of regular polytopes may be found in~\cite{ARP}.

\smallskip
We now turn to chiral polytopes, for which two good references are~\cite{ChiralPolytopes,ChiralPolytopesWithPSL}.

\smallskip
A polytope $\P$ said to be {\em chiral\/} if its automorphism group $G$ has two orbits on its set of flags,
with every two adjacent flags lying in different orbits.
(Another way of viewing this definition is to consider $\P$ as admitting no `reflecting' automorphism
that interchanges a flag with an adjacent flag.)
Here the number of flags of $\P$ is $2|G|$,  and $G$ acts regularly on each of the two orbits.
Again $\P$ is equivelar, with the stabiliser in $G$ of every 2-section of $\P$ inducing the full cyclic group on the corresponding polygon if $\P$ is a polyhedron and the full dihedral group otherwise.


 For a given flag $\Phi$, denote by $F_i$  the $i$-face in $\Phi$ for each $0\leq i\leq n$. For every $j \in \{1,2,\dots,n-1\}$, the chiral polytope $\P$ admits an automorphism $\s_j$ that fixes each $F_i$ with $i\not=j-1,j$ and cyclically permutes
consecutive $j$- and $(j-1)$-faces in the $2$-section $F_{j+1}/F_{j-2}$, 
that is, $\s_j$ takes $\Phi$ to the flag $(\Phi)^{j, j-1}$ which differs from $\Phi$ in precisely its $(j-1)$-and $j$-faces.
This automorphism $\sigma_j$ is the analogue of the abstract rotation $\r_{j-1}\r_{j}$ in the
regular case, for each $j$. These automorphisms $\s_1,\s_2,\dots,\s_{n-1}$ generate $G$, and if $\P$ has Schl\"afli type $\{k_1, k_2, \dots, k_{n-1}\}$,
then the $\s_i$'s satisfy the defining relations for the orientation-preserving subgroup of (index $2$ in) the
string Coxeter group $[k_1, k_2, \cdots, k_{n-1}]$.
Also they satisfy a chiral form of the intersection condition, which is a variant of the one mentioned earlier for regular polytopes.

Chiral polytopes occur in pairs (or {\em enantiomorphic\/} forms), such that each member of the pair is the `mirror image' of the other.
Suppose one of them is $\P$, and has Schl\"afli type $\{k_1, k_2, \cdots, k_{n-1}\}$. Then $G$ is isomorphic
to the quotient of the orientation-preserving subgroup $\Lambda^{\rm o}$ of the string Coxeter group $\Lambda = [k_1, k_2, \cdots, k_{n-1}]$
via some normal subgroup $K$.  By chirality, $K$ is not normal in the full Coxeter group $\Lambda$, but is conjugated by
any orientation-reversing element $c \in \Lambda$ to another normal subgroup $K^c$ which is the kernel of an epimorphism from $\Lambda^{\rm o}$
to the automorphism group $G^c$ of the mirror image $\P^c$ of $\P$.

The automorphism groups of $\P$ and $\P^c$ are isomorphic to each other, but their canonical generating sets
satisfy different defining relations.   In fact, replacing the elements $\s_1$ and $\s_2$ in the canonical generating
tuple $(\s_1,\s_2,\s_3,\dots,\s_{n-1})$ by $\s_1^{-1}$ and $\s_1^{\,2}\s_2$ gives a set of generators for $G$
that satisfy the same defining relations as a canonical generating tuple for $G^c$,
but chirality ensures that there is no automorphism of $G$ that takes $(\s_1,\s_2)$ to $(\s_1^{-1},\s_1^{\,2}\s_2)$
and fixes all the other $\s_j$'s.

Conversely, any finite group $G$ that is generated by $n-1$ elements $\s_1,\s_2,\dots,\s_{n-1}$ which satisfy both
the defining relations for $\Lambda^{\rm o}$ and the chiral form of the intersection condition is the `rotation subgroup' of an
abstract $n$-polytope $\P$ that is either regular or chiral.
Indeed, $\P$ is regular if and only if $G$ admits a group automorphism $\alpha$ of order $2$
that takes $(\s_1,\s_2,\s_3,\dots,\s_{n-1})$ to $(\s_1^{-1},\s_1^{\,2}\s_2,\s_3,\dots,\s_{n-1})$.

\smallskip
We now focus our attention on the rank $3$ case.
Here the generators $\s_1, \s_2$  for $\Aut(\P)$ satisfy the canonical relations
$\,\s_1^{k_1} = \s_2^{k_2}  = (\s_1\s_2)^2 = 1$,
and the chiral form of the intersection condition can be abbreviated to $\lg \s_1\rg \cap \lg \s_2\rg =\{1\}$.
Moreover, $\P$ is chiral if and only if there is no $\r \in \Aut(\P)$ with  $\alpha(\s_1, \s_2) = (\s_1^{-1}, \s_2^{-1})$.

The following proposition is useful for the groups we will deal with in the proof of our main theorem.
It is called the {\em quotient criterion} for chiral $3$-polytopes.

\begin{prop}{\rm \cite[Lemma 3.2]{ChiralPolytopesofPSLextension}}\label{quotient criterion}
$\,$
Let $G$ be a group generated by elements $\s_1, \s_2$ such that $(\s_1\s_2)^2=1$,
and let $\theta\!: G \to H$ be a group homomorphism taking $(\s_1, \s_2)$ to $(\ld_1, \ld_2)$,
such that the restriction of $\theta$ to either $\lg \s_1\rg$ or $\lg \s_2\rg$ is injective.
If $(\ld_{1}, \ld_{2})$ is a canonical generating pair for $H$ as the automorphism group of some chiral $3$-polytope,
then the pair $(\s_1, \s_2)$ satisfies the chiral form of the intersection condition for $G$.
\end{prop}

\subsection{Group theory}

In this section we briefly describe some of the knowledge of group theory we need.
We use standard notation for group theory, as in~\cite{GroupBookss,Huppert,XQ} for example.

Let $G$ be a group.   We define the {\em commutator\/} $[x, y]$ of elements $x$ and $y$ of $G$
by $[x, y]=x^{-1}y^{-1}xy$. 
The group $G':=\lg [x, y] |x, y \in G\rg$ is the {\em commutator subgroup} of $G$.
We let $G_3:=[G', G]=\lg [x, y]| x \in G', y \in G\rg$ and $G_4:=[G_3, G]=\lg [x, y] | x \in G_3, y \in G\rg$.

The following results are elementary and so we give them without proof.

\begin{prop}\label{commutator}
Let $G$ be a group. Then, for any $x, y, z \in G$,
\begin{itemize}
\item [{\rm(1)}] $[xy, z]=[x, z]^y[y, z]$, $[x, yz]=[x, z][x, y]^z$;
\item [{\rm(2)}] $[x, y^{-1}]^y=[x, y]^{-1}$, $[x^{-1}, y]^x=[x, y]^{-1}$;
\item [{\rm(3)}] $[x^{-1}, y^{-1}]^{xy}=[x, y]$;
\item [{\rm(4)}] $[x^{i_1}, y^{i_2}, z^{i_3}]=[x,y,z]^{i_1i_2i_3} ~\quad (\mod G_4)$.
\end{itemize}
\end{prop}

A finite group $G$ is called {\em metacyclic} if it has a cyclic normal subgroup $N$ such that $G/N$ is also cyclic.
Huppert proved the following theorem.

\begin{theorem}\cite[Theorem 11.5]{Huppert}\label{huppert-metacyclic}
Let $G=\lg x\rg \lg y\rg$ be a $p$-group with $p>2$. Then $G$ is metacyclic.
\end{theorem}

A finite $p$-group $G$ is a {\em modular $p$-group} if and only if $HK=KH$, $\forall H, K \leq G$.
Iwaswa proved the following theorem.
\begin{theorem}\cite[Theorem 2.3.1]{1994+Schmidt}\label{modular p group}
Let $G$ be a metacyclic $p$-group with $p>2$. Then $G$ is a modular $p$-group.
\end{theorem}

Let $G$ be a group.
The {\em Frattini subgroup}, denoted by $\Phi(G)$, is the intersection of all maximal subgroups of $G$. Obviously, $\Phi(G)$ is a characteristic subgroup\footnote{Recall that a characteristic subgroup $H$ of a group $G$ is a subgroup that is fixed by the automorphism group of $G$.} of $G$. The following theorem is the well-known Burnside Basis Theorem.

\begin{theorem}{\rm ~\cite[Theorem 1.12]{GroupBookss}}\label{burnside}
Let $G$ be a $p$-group and $|G: \Phi(G)| = p^d$.
\begin{itemize}
\item [(1)] $G/\Phi(G) \cong \mz_p^d$. Moreover, if $N \lhd G$ and $G/N$ is elementary abelian, then $\Phi(G) \leq N$.

\item [(2)] Every minimal generating set of $G$ contains exactly $d$ elements.
\end{itemize}
\end{theorem}

The unique cardinality of all  minimal generating sets of a $p$-group $P$ is called the {\em rank} of $P$,  and denoted by $d(P)$.

We will also need Sylow's third theorem so we state it here for clarity.
\begin{theorem}[3rd Sylow's theorem]\label{sylow}
Let $G$ be a group of order $p^nm$ with $(p,m)=1$. Let $n_p(G)$ be the number of Sylow $p$-subgroup of $G$. Then
\begin{enumerate}
    \item $n_p(G)\; | \;m$;
    \item $n_p(G) \equiv 1\; \mod p$;
    \item $n_p(G) = [G:N_G(P)]$ where $P$ is a Sylow $p$-subgroup of $G$.
\end{enumerate}\end{theorem}

Let $c(P)$ denote the nilpotency class of a group $P$. From the classification of finite groups of order $p^3$ and $p^4$ (see~\cite[Chapter 4]{Hall} and~\cite[Theorem 12.6]{Huppert}) we have the following propositions. 

\begin{prop}\label{group of order $p^3$}
Let $P$ be a non-abelian group of order $p^3$ for odd prime $p>2$. Suppose $P$ is non-metacyclic.
Then $P \cong M_{p}(1,1,1)$, and $c(P)=2, P'=Z(P) \cong \mathbb{Z}_{p}$, where
$$M_{p}(1,1,1)=\lg a, b, c \ | \ a^p, b^p, c^p, c=[a,b], [a,c], [b, c]\rg\cong (\mathbb{Z}_p \times \mathbb{Z}_p) \rtimes \mathbb{Z}_p.$$
\end{prop}

\begin{prop}\label{group of order $p^4$}
Let $P$ be a non-abelian group of order $p^4$ for odd prime $p>3$. Suppose $P$ is non-metacyclic and $d(P)=2$.
Then $c(P)=2$ or $3$. Furthermore, the following hold:
\begin{itemize}
\item [\rm (1)] if $c(P)=2$, then $P' \cong \mathbb{Z}_{p}$, and  $\Omega_1(P)=\lg x \in P \ | \ x^p=1\rg <P$;
\item [\rm (2)] if $c(P)=3$, then $P' \cong \mathbb{Z}_{p} \times \mathbb{Z}_{p}$, and one of the following holds:
\begin{itemize}
\item [\rm(a)] $\exp(P)=p$, $Z(P)=P_3 \cong \mathbb{Z}_{p}$;

\item [\rm (b)] $\exp(P)=p^2$, $Z(P)=P_3=\mho_1(P) \cong \mathbb{Z}_{p}$. Moreover, $\mho_1(P)=\lg x^p | x \in P, o(x)=p^2\rg$.

\end{itemize}

\end{itemize}

\end{prop}

Finally, we state a proposition about binomial coefficients that will be used in the proof of one of our theorems.

\begin{prop}\cite[Chapter 1]{MA}\label{enumeration}
Let $n, m, k \in \mathbb{N}$ and let $\tbinom{n}{m}$ denote the binomial coefficient indexed by $n$ and $m$. Then
\begin{itemize}
\item [\rm(1)] $\tbinom{n}{m} =\tbinom{n}{n-m}$;
\item [\rm(2)] $\tbinom{n+2}{m+1} =\tbinom{n+1}{m+1}+\tbinom{n+1}{m}$;
\item [\rm(3)] $\tbinom{n}{m}\tbinom{m}{k}=\tbinom{n}{k}\tbinom{n-k}{m-k}$;
\item [\rm(4)] $\sum\limits_{k=0}^m(-1)^{k}\tbinom{n}{k}=(-1)^m\tbinom{n-1}{m}$;
\end{itemize}
Furthermore, if $m <0$, then $\tbinom{n}{m}=0$.
\end{prop}

\section{A few technical lemmas}\label{sec3}

In this section, we give lemmas about finite $p$-groups of maximal class that are useful in this paper. 
We assume in what follows that $p$ is an odd prime and $m$ is a strictly positive integer. A group $P$ of order $p^m$ and nilpotency class $m-1$ is said to be of {\em maximal class} if $m \geq 3$. 
The basic material about these groups can be found in Blackburn~\cite[pages 83--84]{Blackburn} or Huppert~\cite[Chapter 3]{Huppert}.

Let $P$ be a group of maximal class and order $p^m$. If $P$ has an abelian maximal subgroup $A$, and any $x \in P \backslash   A$ has order $p^2$, 
then, by \cite[Chapter 8, Example 8.3.4]{XQ}), 

$$P=\lg s_1, s_2, \cdots, s_r, s_{r+1}, \cdots, s_{p-1}, \b \rg\ \mbox{with}\ \mbox{the}\ \mbox{following}\ \mbox{relations}$$
$$ s_1^{p^e}, s_2^{p^e}, \cdots, s_r^{p^e}, s_{r+1}^{p^{e-1}}, \cdots s_{p-1}^{p^{e-1}}, $$
$$\b^{p^2}, \b^p=s_r^{p^{e-1}},$$
$$[s_i, s_j],  \ {\rm for} \ 1 \leq i, j \leq p-1,\ \mbox{and}\ s_{k+1}=[s_k, \b],  \ {\rm for} \ 1 \leq k \leq p-2,$$
$$[s_{p-1}, \b]=s_1^{-\tbinom{p}{1}}s_2^{-\tbinom{p}{2}}\cdots s_{p-2}^{-\tbinom{p}{p-2}}s_{p-1}^{-\tbinom{p}{p-1}},$$
where 
$$1 \leq e,1 \leq r \leq p-1,\ m=(p-1)(e-1)+r+1 \geq 3.$$

\noindent Noting that  $A = \lg s_1, s_2, \cdots, s_{p-1}\rg$, the group $A$ is abelian and $A \unlhd P$. 
Moreover, we also have $(\b\sigma_1)^p=\b^p=s_{r}^{p^{e-1}}$ and $o(\b\sigma_1)=p^2$.
In fact, any $x \in P\backslash  A$ has order $p^2$.
With the relations given above, it is clear that $P=\lg s_1, \b\rg$. 

Define
\begin{center}
\begin{tabular}{lll}
$\tau:$ \qquad & $s_1 \mapsto s_1^{-1}$ \quad &  and \quad $\beta \mapsto \beta$\\
$\sigma:$ \qquad & $s_1 \mapsto s_1$ \quad &  and \quad $\beta \mapsto \beta^{-1}$. 
\end{tabular}
\end{center}

For convenience, let $s_{k+1}:=[s_{k}, \b]$ for $k = p-1$ and $p$. Then
$$s_p=[s_{p-1}, \b]=\prod\limits_{i=1}^{p-1}s_i^{-\tbinom{p}{i}},$$
$$s_{p+1}=[s_p,\b]=s_p^{-1}s_p^{\b}=\prod\limits_{i=1}^{p-1}(s_i^{-1})^{-\tbinom{p}{i}}\cdot\prod\limits_{i=1}^{p-1}(s_i^{\beta})^{-\tbinom{p}{i}}$$
$$=\prod\limits_{i=1}^{p-1}(s_i^{-1}s_i^{\beta})^{-\tbinom{p}{i}}=\prod\limits_{i=1}^{p-1}[s_i, \beta]^{-\tbinom{p}{i}}=\prod\limits_{i=1}^{p-1}s_{i+1}^{-\tbinom{p}{i}}.$$

We prove the following lemma.

\begin{lem}
$\tau \notin \Aut(P)$.
\end{lem}

\demo Suppose $\tau \in \Aut(P)$. Since $A$ is abelian and $s_{i+1}=[s_i, \b]$ for $1 \leq i \leq p-2$, by Proposition~\ref{commutator}, we have
$$s_{i+1}^{\tau}=[s_i, \b]^{\tau}=[s_i^{-1}, \b]=(([s_i, \b])^{-1})^{s_i^{-1}}=(s_{i+1}^{-1})^{s_i^{-1}}=s_{i+1}^{-1}.$$
In particular, $(s_r^{p^{e-1}})^{\tau}=(s_r^{\tau})^{p^{e-1}}=(s_r^{p^{e-1}})^{-1}$.
On the other hand, since $\beta^{p}=s_{r}^{p^{e-1}}$, $\beta^p=(\beta^{p})^{\tau}=(s_{r}^{p^{e-1}})^{\tau}=\b^{-p}$. 
This implies that $\b^{2p}=1$. But this is impossible because $o(\b)=p^2$. \qed

\begin{lem}\label{lemma1}\label{technical lemma3.2}
For any $1 \leq k \leq p$, we have
\begin{itemize}
\item [{\rm(1)}]$s_k^{\b^{1-k}}=\prod\limits_{i=k}^{p+1}s_i^{\tbinom{p+1-k}{i-k}}$;
\item [{\rm(2)}]$s_{k}^{\sigma}=\left\{
\begin{array}{ll}
s_k^{\b^{1-k}},   & k  \mbox{ is odd,} \\
(s_k^{\b^{1-k}})^{-1},  & k   \mbox{ is even.}
\end{array}
\right.$

\end{itemize}
\end{lem}

\demo
(1) Since $s_{k+1}=[s_k, \b]$, we have $s_k^{\b}=s_ks_{k+1}$.
Since $o(\b)=p^2, \b^{p} \in A$ and $k-1 \leq p$, we have $p-k+1 \geq 0$.
It follows that
$$s_k^{\b^{1-k}}=s_k^{\b^{p+1-k}}=(s_ks_{k+1})^{\b^{p-k}}=(s_ks_{k+1}^2s_{k+2})^{\b^{p-k-1}}=\cdots $$
$$=s_k^{\tbinom{p+1-k}{0}}s_{k+1}^{\tbinom{p+1-k}{1}}s_{k+2}^{\tbinom{p+1-k}{2}}\cdots s_{p}^{\tbinom{p+1-k}{p-k}}s_{p+1}^{\tbinom{p+1-k}{p+1-k}}.$$

\f (2) The proof uses an induction on $k$. For $k=1$, we have $s_1^{\sigma}=s_1$. For $k=2$, by Proposition~\ref{commutator}, we have $$s_2^{\sigma}
=[s_1, \b]^{\sigma}=[s_1, \b^{-1}]=([s_1, \b]^{-1})^{\b^{-1}}=(s_2^{-1})^{\b^{1-2}}.$$

\noindent Hence this lemma is true for $k=1$ and $2$. We now use induction on $k$ and Proposition~\ref{commutator}, splitting the analysis in two cases, namely the case where $k$ is odd and the case where $k$ is even.
 
If $k$ is odd, then $$s_{k}^{\sigma}=[s_{k-1}^{\sigma}, \b^{\sigma}]=[(s_{k-1}^{\b^{2-k}})^{-1}, \b^{-1}]=[s_{k-1}^{-1}, \b^{-1}]^{\b^{2-k}}=([s_{k-1}, \b]^{\b^{-1}s_{k-1}^{-1}})^{\b^{2-k}}.$$
Recall that $A=\lg s_1, s_2, \cdots, s_{p-1}\rg$ is an abelian maximal subgroup of $G$ and $A \unlhd G$.
It follows that $s_k \in A$ and $s_k^{\b} \in A$ for any $k \geq 1$. Thus,
$$s_{k}^{\sigma}=([s_{k-1}, \b]^{\b^{-1}s_{k-1}^{-1}})^{\b^{2-k}}=s_k^{(s_{k-1}^{-1})^{\b}\b^{1-k}}=s_k^{\b^{1-k}}.$$

If $k$ is even, then
$$s_{k}^{\sigma}=[s_{k-1}^{\sigma}, \b^{\sigma}]=[s_{k-1}^{\b^{2-k}}, \b^{-1}]=[s_{k-1}, \b^{-1}]^{\b^{2-k}}=(([s_{k-1}, \b]^{-1})^{\b^{-1}})^{\b^{2-k}}=(s_k^{\b^{1-k}})^{-1}.$$
\qed

\begin{lem}~\label{reven}
If $m$ is odd , then $r$ is even and $\sigma \in \Aut(P)$. Moreover, $o(\sigma)=2$.
\end{lem}

\demo Since $m=(p-1)(e-1)+r+1$ and $p$ and $m$ are odd, $r$ must be even. 
The fact that $o(\sigma) = 2$ follows from the definition of $\sigma$.
The group $P = \lg s_1^{\s}, \b^{\s}\rg=\lg s_1, \b^{-1}\rg$. Thus, we only need to show that the generating set $\{s_1^{\s}, s_2^{\s}, \cdots, s_{p-1}^{\s}, \b^{\s}\}$
of $P$ satisfies
the same relations as $\{s_1, s_2, \cdots, s_{p-1}, \b\}$ in order to conclude that $\sigma \in \Aut(P)$.
Recall that $A=\lg s_1, s_2, \cdots, s_{p-1}\rg$ is abelian and $A \unlhd P$. By Lemma~\ref{technical lemma3.2},
$$s_{i}^{\sigma}=\left\{
\begin{array}{ll}
s_i^{\b^{1-i}},   & \mbox{ for $i$ odd,} \\
(s_i^{\b^{1-i}})^{-1},  & \mbox{ for $i$ even.}
\end{array}
\right.$$
Now $o(s_i)=o(s_i^{\s})$, $[s_i^{\s}, s_j^{\s}]=1$ for $1 \leq i, j \leq p-1$.
Moreover, it is easy to see that $s_{k+1}^{\sigma}=[s_k^{\s}, \b^{\s}]$ for $1\leq k \leq p-2$. Thus we only need to show that $(\b^{p})^{\sigma}=(s_r^{p^{e-1}})^{\sigma}$, and
$$[s_{p-1}^{\sigma}, \b^{\sigma}]=(s_1^{\sigma})^{-\tbinom{p}{1}}(s_2^{\sigma})^{-\tbinom{p}{2}}(s_3^{\sigma})^{-\tbinom{p}{3}}\cdots (s_{p-1}^{\sigma})^{-\tbinom{p}{p-1}},$$
that is, $s_p^{\sigma}=\prod\limits_{i=1}^{p-1}(s_i^{\sigma})^{-\tbinom{p}{i}}$.
Since $r$ is even, 
we have
$$(s_r^{p^{e-1}})^{\sigma}=(s_r^{\sigma})^{p^{e-1}}=((s_r^{\b^{1-r}})^{-1})^{p^{e-1}}
=((s_r^{p^{e-1}})^{\b^{1-r}})^{-1}=((\b^p)^{\b^{1-r}})^{-1}=\b^{-p}=(\b^{p})^{\sigma}.$$
Next, by Lemma~3.2, we have
$$s_1^{\sigma}=s_1, s_2^{\sigma}=\prod\limits_{i=2}^{p+1}s_i^{-\tbinom{p-1}{i-2}},
s_3^{\sigma}=\prod\limits_{i=3}^{p+1}s_i^{\tbinom{p-2}{i-3}}, \cdots, s_{p-1}^{\sigma}=
\prod\limits_{i=p-1}^{p+1}s_i^{-\tbinom{2}{i-(p-1)}}.$$

\f Then $\prod\limits_{i=1}^{p-1}(s_i^{\sigma})^{-\tbinom{p}{i}}
=s_1^{u_1}s_2^{u_2}\cdots s_{p+1}^{u_{p+1}}$, where $u_i$ comes from Table~\ref{tableU}.
\begin{table}[]
    \centering
\setlength{\arraycolsep}{5pt}
\begingroup
\renewcommand*{\arraystretch}{1.8}
$$
 \begin{bmatrix}

   \begin{array}{c | c|cccccc}
*            &(s_1^{\sigma})^{-\tbinom{p}{1}}   &(s_2^{\sigma})^{-\tbinom{p}{2}}   &(s_3^{\sigma})^{-\tbinom{p}{3}}   &(s_4^{\sigma})^{-\tbinom{p}{4}}       &\cdots &\cdots &(s_{p-1}^{\sigma})^{-\tbinom{p}{p-1}}\\ \hline
u_1       &-\tbinom{p}{1}      & 0                                                & 0                                                     &0 &\cdots &\cdots         &0                               \\ \hline
u_2       &0&\tbinom{p-1}{0}\tbinom{p}{2}      & 0                                                & 0                                                     &\cdots &\cdots &0                                      \\
u_3        &0&\tbinom{p-1}{1}\tbinom{p}{2}      & -\tbinom{p-2}{0}\tbinom{p}{3}      & 0                                                    &\cdots &\cdots &0                                       \\
u_4        &0&\tbinom{p-1}{2}\tbinom{p}{2}      & -\tbinom{p-2}{1}\tbinom{p}{3}      & \tbinom{p-3}{0}\tbinom{p}{4}        &\cdots &\cdots &0                                        \\
u_5      &0&\tbinom{p-1}{3}\tbinom{p}{2}      & -\tbinom{p-2}{2}\tbinom{p}{3}      & \tbinom{p-3}{1}\tbinom{p}{4}        &\cdots &\cdots &0                                     \\
\vdots    &\vdots&\vdots                                          & \vdots                                         & \vdots                                            &\cdots &\cdots &  \vdots                             \\
\vdots    &\vdots&\vdots                                          & \vdots                                         & \vdots                                            &\cdots &\cdots &   \vdots                                \\
u_{p-1}  &0&\tbinom{p-1}{p-3}\tbinom{p}{2}   & -\tbinom{p-2}{p-4}\tbinom{p}{3}   & \tbinom{p-3}{p-5}\tbinom{p}{4}     &\cdots &\cdots &\tbinom{2}{0}\tbinom{p}{p-1}   \\ \hline
u_p     &0&\tbinom{p-1}{p-2}\tbinom{p}{2}   & -\tbinom{p-2}{p-3}\tbinom{p}{3}   & \tbinom{p-3}{p-4}\tbinom{p}{4}     &\cdots &\cdots &\tbinom{2}{1}\tbinom{p}{p-1} \\ \hline
u_{p+1}&0&\tbinom{p-1}{p-1}\tbinom{p}{2}   & -\tbinom{p-2}{p-2}\tbinom{p}{3}   & \tbinom{p-3}{p-3}\tbinom{p}{4}     &\cdots &\cdots &\tbinom{2}{2}\tbinom{p}{p-1}

  \end{array}

  \end{bmatrix}
$$
\endgroup
    \caption{The values of $u_i$.}
    \label{tableU}
\end{table}
This gives $u_1=-\tbinom{p}{1}$. By Proposition~\ref{enumeration}(1), we have
$$u_{p+1}=\sum\limits_{j=2}^{p-1}(-1)^{j}\tbinom{p}{j}=\tbinom{p}{p-1}+
\sum\limits_{j=2}^{\frac{p-3}{2}}\left[(-1)^{j}\tbinom{p}{j}+(-1)^{p-j}\tbinom{p}{p-j}\right]=\tbinom{p}{p-1}=p.$$
\f For $2 \leq i \leq p-1$, we have $u_i=\sum\limits_{j=2}^{i}(-1)^j\tbinom{p+1-j}{i-j}\tbinom{p}{j}$.
For $i=p$, we have $$u_p-\tbinom{1}{0}\tbinom{p}{p}=u_p-1=
\sum\limits_{j=2}^{p}(-1)^j\tbinom{p+1-j}{p-j}\tbinom{p}{j}.$$
By Proposition~\ref{enumeration}$(2)$$(3)$, we have
$$\tbinom{p+1-j}{i-j}\tbinom{p}{j}=\left[\tbinom{p-j}{i-j}+\tbinom{p-j}{i-j-1}\right]\tbinom{p}{j}
=\tbinom{p-j}{i-j}\tbinom{p}{j}+\tbinom{p-j}{i-j-1}\tbinom{p}{j}=
\tbinom{p}{i}\tbinom{i}{j}+\tbinom{p}{i-1}\tbinom{i-1}{j}.$$
Together with Proposition~\ref{enumeration}(4), for $2 \leq i \leq p$, we have
$$\sum\limits_{j=2}^{i}(-1)^j\left[\tbinom{p}{i}\tbinom{i}{j}+\tbinom{p}{i-1}\tbinom{i-1}{j}\right]$$
$$=\tbinom{p}{i}\left[i-1+\sum\limits_{j=0}^{i}(-1)^j\tbinom{i}{j}\right]+
\tbinom{p}{i-1}\left[i-2+\sum\limits_{j=0}^{i}(-1)^j\tbinom{i-1}{j}\right]$$
$$=\tbinom{p}{i}\left[i-1+(-1)^i\tbinom{i-1}{i}\right]+
\tbinom{p}{i-1}\left[i-2+(-1)^i\tbinom{i-2}{i}\right]$$
$$=\tbinom{p}{i}(i-1)+\tbinom{p}{i-1}(i-2)$$
$$=-\tbinom{p}{i}+(p-i+1)\tbinom{p}{i-1}+(i-2)\tbinom{p}{i-1}$$
$$=-\tbinom{p}{i}
+(p-1)\tbinom{p}{i-1},$$
and hence
$$u_p=-\tbinom{p}{p}
+(p-1)\tbinom{p}{p-1}+1=(p-1)\tbinom{p}{p-1},$$
$$u_i=-\tbinom{p}{i}
+(p-1)\tbinom{p}{i-1}, \quad 2 \leq i \leq p-1.$$
By Lemma~3.2 and the fact that $s_{p+1}=[s_p,\b]$, we have $s_p^{\sigma}=s_p^{\b^{1-p}}=s_p^{\b}=s_ps_{p+1}.$
The result is now immediate since $$s_1^{u_1}s_2^{u_2}\cdots s_p^{u_p}s_{p+1}^{u_{p+1}}$$
$$=s_1^{-\tbinom{p}{1}}\left[\prod\limits_{i=2}^{p-1}s_i^{-\tbinom{p}{i}
+(p-1)\tbinom{p}{i-1}}\right]s_{p}^{(p-1)\tbinom{p}{p-1}}s_{p+1}^{p}=\left[\prod\limits_{i=1}^{p-1}s_i^{-\tbinom{p}{i}}\right]\left[(\prod\limits_{i=2}^{p}s_i^{-\tbinom{p}{i-1}})^{1-p}\right]s_{p+1}^p$$
$$=s_ps_{p+1}^{1-p}s_{p+1}^{p}=s_{p}s_{p+1}=s_p^{\sigma}.$$ 
\qed

\section{Proof of Theorem~\ref{maintheorem1}}\label{sec4}

(1) Let $n_p(G)$ be the number of Sylow $p$-subgrougs of $G$ and let $P$ be a Sylow $p$-subgroup of $G$.
By Theorem~\ref{sylow} (1),
$n_p(G) \mid 2$ and $n_p(G) \equiv  1\ \mod p$. It follows that $n_p(G)=1$.
Then $P \unlhd G$. 

Since $| \lg \sigma_1\sigma_2 \rg|=2$ and $|P|=p^m$, $P \cap \lg \sigma_1\sigma_2\rg=\{1\}$.
It follows that $|P \lg \sigma_1\sigma_2 \rg|=2p^m=|G|$ and thus $G=P\lg \sigma_1\sigma_2\rg$.
In face, since $P \unlhd G$, $G=P \rtimes \lg \sigma_1\sigma_2\rg$, finishing the proof of (1).

(2) Up to duality, the Schl\"afli type $\{k_1, k_2\}$  must be one of $\{p^{l_1}, p^{l_2}\}, \{p^{l_1}, 2p^{l_2}\}$ or $\{2p^{l_1}, 2p^{l_2}\}$. 
Consider the quotient group $\overline{G}=G/P=\lg \overline{\sigma_1}, \overline{\sigma_2}\rg$, then $|\overline{G}|=2$.
If $\{k_1, k_2\}=\{p^{l_1}, p^{l_2}\}$, then $\sigma_1, \sigma_2 \in P$, and hence $\overline{\sigma_1}=\overline{\sigma_2}=\overline{1}$ in $\overline{G}$,
which is impossible because $|\overline{G}|=2$. If $\{k_1, k_2\}=\{2p^{l_1}, 2p^{l_2}\}$, then $\overline{\sigma_1} \ne \overline{1}$, $\overline{\sigma_2} \ne \overline{1}$ and
$(\overline{\sigma_1})^2 =(\overline{\sigma_2})^2=\overline{1}$. Since $(\sigma_1\sigma_2)^2=1$ in $G$, $\overline{\sigma_1} \overline{\sigma_2} \ne \overline{1}$ and 
$(\overline{\sigma_1}\overline{\sigma_2})^2=\overline{1}$ in $\overline{G}$. Thus $\overline{G}=\lg \overline{\sigma_1}, \overline{\sigma_2}\rg \cong \mathbb{Z}_2 \times \mathbb{Z}_2$, which is impossible because $|\overline{G}|=2$. Therefore,  $\{k_1, k_2\}=\{p^{l_1}, 2p^{l_2}\}$, finishing the proof of (2).

(3) By Theorem~\ref{maintheorem1}$(2)$, we have  $\{k_1, k_2\}=\{p^{l_1}, 2p^{l_2}\}$. Let $H=\lg \sigma_1, \sigma_2^2\rg$. Then $H \leq P$. Obviously, $\sigma_1 \in H \leq N_{G}(H)$. Moreover,
since $(\sigma_1\sigma_2)^2=1$, 
$$\sigma_1^{\sigma_2}=\sigma_2^{-1}\sigma_1\sigma_2=\sigma_2^{-1}\sigma_2^{-1}\sigma_1^{-1}=\sigma_2^{-2}\sigma^{-1} \in H, (\sigma_2^2)^{\sigma_2}=(\sigma_2)^2 \in H.$$
Then $\sigma_2 \in N_{G}(H)$. It follows that $H\unlhd G$ and $H\lg \sigma_1\sigma_2\rg=\lg \sigma_1, \sigma_2^2, \sigma_1\sigma_2\rg= G$.
On the other hand, it is easy to see that $|G/H|=2$. Thus, $P=H$ and $d(P)=d(H) \leq 2$. 

If $d(P)=1$, then $P$ is cyclic. Since $p \mid k_1$ and $p \mid k_2$, we have $\sigma_1, \sigma_2^2 \in P$, and so $\lg \sigma_1 \rg \leq \lg \sigma_2^2\rg$ or 
$\lg \sigma_1 \rg \geq \lg \sigma_2^2\rg$.
It follows that $|\lg \sigma_1 \rg \cap \lg \sigma_2^2\rg| >1$, which is impossible
 because $\lg \sigma_1 \rg \cap \lg \sigma_2 \rg=\{1\}$.
Thus, $d(P)=2$. 

Now, suppose that $P$ is abelian. Then $[\sigma_1, \sigma_2^2]=1$ in $P$ and hence $P=\lg \sigma_1, \sigma_2^2\rg=\lg \sigma_1\rg \lg \sigma_2^2\rg$. 
Since $\lg \sigma_1 \rg \cap \lg \sigma_2\rg=\{1\}$, we have 
$$p^m=|P|=|\lg \sigma_1, \sigma_2^2 \rg|=|\lg \sigma_1 \rg \lg \sigma_2^2\rg|=\frac{|\lg \sigma_1\rg| |\lg \sigma_2^2\rg|}{|\lg \sigma_1\rg \cap \lg \sigma_2^2\rg|}=p^{l_1}p^{l_2}=p^{l_1+l_2}.$$
It means that $\mathcal{P}$ is tight, that is $|G|=2p^{l_1+l_2}$.  Moreover, we claim that $G=\mathcal{G}$, where 
$$\mathcal{G}=\lg \sigma_1, \sigma_2 \ | \ \sigma_1^{p^{l_1}} , \sigma_2^{2p^{l_2}}, (\sigma_1\sigma_2)^2, [\sigma_1, \sigma_2^2]\rg.$$
Obviously, $(G, \{\sigma_1, \sigma_2\})$ satisfies all of the defining relations of $(\mathcal{G}, \{\sigma_1, \sigma_2\})$, and so $G$ is the quotient group of $\mathcal{G}$. On the other hand,
since $[\sigma_1, \sigma_2^2]=1$, $\lg \sigma_2^2\rg \unlhd \mathcal{G}$. Consider the quotient group 
$$\overline{\mathcal{G}}=\mathcal{G}/\lg \sigma_2\rg=\lg \overline{\sigma_1}, \overline{\sigma_2}  \ | \ \overline{\sigma_1}^{p^{l_1}}, \overline{\sigma_2}^{2}, (\overline{\sigma_1}\overline{\sigma_2})^2\rg.$$
It is easy to see that $\overline{\mathcal{G}} \leq D_{2p^{l_1}}$ where $D_{2p^{l_1}}$ is the dihedral group of order $2p^{l_1}$. Then
$|\mathcal{G}| =  |\overline{\mathcal{G}}|\cdot |\lg \sigma_2^2 \rg|\leq 2{p^{l_1}}p^{l_2}=2p^{l_1+l_2}$. So $|\mathcal{P}|=2p^{l_1+l_2}=|G|$, and hence $G=\mathcal{G}$.
However, consider the map $\alpha: \sigma_1 \mapsto \sigma_1^{-1}, \sigma_2 \mapsto \sigma_2^{-1}$. It is easy to check that 
$$(\sigma_1^{-1})^{p^{l_1}}=(\sigma_{2}^{-1})^{2p^{l_2}}=(\sigma_1^{-1}\sigma_2^{-1})^2=[\sigma_1^{-1}, ((\sigma_2^{-1})^2)]=1.$$
Then $\alpha \in \Aut(G)$, which is impossible because $\mathcal{P}$ is chiral. Thus, $P$ is not abelian and $m \geq 3$.
This finishes the proof of (3).

4) Suppose first that $\mathcal{P}$ is tight, that is $|G| = k_1k_2 = 2p^m$. Then $|P| = o(\sigma_1)o(\sigma_2^2)$ and since $\lg \sigma_1\rg \cap \lg \sigma_2\rg = \{1\}$, we have $P = \lg \sigma_1\rg \lg \sigma_2\rg$. Theorem~\ref{huppert-metacyclic} implies that $P$ is metacyclic group.

Now suppose that $P = \lg \sigma_1, \sigma_2^2\rg$ is metacyclic. By Theorem~\ref{modular p group}, we have $P=\lg \sigma_1\rg\lg \sigma_2^2\rg$.
Note that $\lg \sigma_1\rg \cap \lg \sigma_2^2\rg=\{1\}$, so $|P|=|\lg \sigma_1\rg|| \lg \sigma_2^2\rg|$. 
It follows that $|G|=2\cdot o(\sigma_1) \cdot o(\sigma_2^2)=o(\sigma_1)\cdot o(\sigma_2)=k_1k_2$ and hence $\mathcal{P}$ is tight.

\section{Proof of Theorem~\ref{maintheorem2}}\label{sec5}
For $p=3$ or $5$, we may refer to the lists of all chiral polytopes with up to $4000$ flags computed by Conder~\cite{atlas1}. 
So, in the following discussion, we always assume $p \geq 7$, and $P \in \Syl_p(G)$.

(1) Suppose $|G|=2p^3$ and $\mathcal{P}$ is a non-tight chiral polyhedron. By Theorem~\ref{maintheorem1}(2),  we know that $\{k_1, k_2\}=\{p, 2p\}$.
Then $o(\sigma_1)=o(\sigma_2^2)=p$, which means that $P=\lg \sigma_1, \sigma_2^2\rg$ is generated by two elements of order $p$.
By Proposition~\ref{group of order $p^3$}, $\lg [\sigma_1, \sigma_2^2]\rg=P'=Z(P) \cong \mathbb{Z}_p$ and hence 
$[\sigma_1, \sigma_2^2]^p=[[\sigma_1, \sigma_2^2], \sigma_1]=[[\sigma_1, \sigma_2^2], \sigma_2^2]=1$ in $P$.
Knowing how to represent these two elements of $P$ is enough to define the entire group $G$.
We claim that $G=\mathcal{G}_1$, where 
$$\mathcal{G}_1=\lg \sigma_1, \sigma_2  \ | \ \sigma_1^p, \sigma_2^{2p}, (\sigma_1\sigma_2)^2, [\sigma_1, \sigma_2^2]^p, [[\sigma_1, \sigma_2^2], \sigma_1], [[\sigma_1, \sigma_2^2], \sigma_2^2] \rg.$$
Obviously, $(G, \{\sigma_1, \sigma_2\})$ satisfies all of the defining relations of $(\mathcal{G}_1, \{\sigma_1, \sigma_2\})$, and hence $G$ is a quotient group of $\mathcal{G}_1$. 
On the other hand, since $(\sigma_1\sigma_2)^2=1$ in $\mathcal{G}_1$, we have $\sigma_1^{\sigma_2}=\sigma_2^{-2}\sigma_1^{-1} \in \lg \sigma_1, \sigma_2^2\rg$ and hence $\lg \sigma_1, \sigma_2^2\rg \unlhd \mathcal{G}_1$. Consider the quotient group 
$\overline{\mathcal{G}_1}=\mathcal{G}_1/\lg \sigma_1, \sigma_2^2\rg=\lg \overline{\sigma_2} \ | \  \overline{\sigma_2}^2 \rg$. Then $|\overline{\mathcal{G}_1}| \leq 2$.
Note that the subgroup $\lg \sigma_1, \sigma_2^2\rg$ of $\mathcal{G}_1$ satisfies all of the defining relations of the subgroup $P$ of $G$, so 
the subgroup $\lg \sigma_1, \sigma_2^2\rg$ of $\mathcal{G}_1$ has order at most $|P|=p^3$. Hence $|\mathcal{G}_1| \leq 2p^3=|G|$. 
Thus, $G=\mathcal{G}_1$.

Consider the map $\alpha: \sigma_1 \mapsto \sigma_1^{-1}, \sigma_2 \mapsto \sigma_2^{-1}$. It is easy to check that 
$$(\sigma_1^{-1})^p=(\sigma_{2}^{-1})^{2p}=(\sigma_1^{-1}\sigma_2^{-1})^2=[\sigma_1^{-1}, ((\sigma_2^{-1})^2)]^p=[\sigma_1^{-1}, (\sigma_2^{-1})^2, \sigma_1^{-1}]=[\sigma_1^{-1}, (\sigma_2^{-1})^2, (\sigma_2^{-1})^2]=1$$
meaning $\alpha \in \Aut(G)$, which contradicts the hypothesis that $\mathcal{P}$ is chiral. Thus, $\mathcal{P}$ must be tight.

(2) Suppose $|G|=2p^4$ and $\mathcal{P}$ is a non-tight chiral polyhedron. By Theorem~\ref{maintheorem1}(4), $P$ is non-metacyclic. 
Moreover, by Theorem~\ref{maintheorem1}(2) and Proposition~\ref{group of order $p^3$}, we have $\exp(P')=p$, $c(P)=2$ or $3$ and $\{k_1, k_2\} \in \{\{p, 2p\}, \{p, 2p^2\}, \{p^2,2p\}\}$, meaning we have six cases to consider, that is three for each possible value of $c(P)$. We show for each of these six cases that they cannot happen.

\noindent {\bf Case 1:} $c(P)=2$. 

\noindent{\bf Subcase 1.1:} $\{k_1, k_2\}=\{p, 2p\}$. 
Then $o(\sigma_1)=o(\sigma_2^2)=p$, and so 
$$P=\lg \sigma_1, \sigma_2^2\rg \leq \Omega_1(P)=\lg x \in P  \ | \ x^p=1\rg.$$
But this is impossible by Proposition~\ref{group of order $p^4$}(1).

 \noindent{\bf Subcase 1.2:} $\{k_1, k_2\}=\{p, 2p^2\}$. Then $o(\sigma_1)=p$ and $o(\sigma_2^2)=p^2$. 
 Moreover, 

 $$[\sigma_1, \sigma_2^2]^p=1, [\sigma_1, \sigma_2^2, \sigma_1]=1, [\sigma_1, \sigma_2^2, \sigma_2^2]=1.$$
 This is enough to define the entire group $G$. We claim that $G=\mathcal{G}_2$, where
$$\mathcal{G}_2= \lg \sigma_1, \sigma_2  \ | \ \sigma_1^p, \sigma_2^{2p^2}, (\sigma_1\sigma_2)^2,
 [\sigma_1, \sigma_2^2]^p, [\sigma_1, \sigma_2^2, \sigma_1], [\sigma_1, \sigma_2^2, \sigma_2^2] \rg.$$
Obviously, $(G, \{\sigma_1, \sigma_2\})$ satisfies all of the defining relations of $(\mathcal{G}_2, \{\sigma_1, \sigma_2\})$, and hence $G$ is a quotient group of $\mathcal{G}_2$. 
Note that  $\lg \sigma_1, \sigma_2^2\rg \unlhd \mathcal{G}_2$. Consider the quotient group 
$\overline{\mathcal{G}_2}=\mathcal{G}_2/\lg \sigma_1, \sigma_2^2\rg=\lg \overline{\sigma_2} \ | \  (\overline{\sigma_2})^2 \rg$. Then $|\overline{\mathcal{G}_2}| \leq 2$.
For the subgroup group $\lg \sigma_1, \sigma_2^2\rg$ of $\mathcal{G}_2$, since $[\sigma_1, \sigma_2^2]^p=[\sigma_1, \sigma_2^2, \sigma_1]=[\sigma_1, \sigma_2^2, \sigma_2^2]=1$, $\lg \sigma_1, \sigma_2^2\rg'=\lg [\sigma_1, \sigma_2^2]\rg$. Then
$$|\lg \sigma_1, \sigma_2^2\rg| =|\lg \sigma_1, \sigma_2^2\rg / \lg \sigma_1, \sigma_2^2\rg'| \cdot |\lg \sigma_1, \sigma_2^2\rg'| \leq p\cdot p^2 \cdot p=p^4.$$
It follows that
 $|\mathcal{G}_2|=|\overline{\mathcal{G}_2}|\cdot | \lg \sigma_1, \sigma_2^2\rg| \leq 2p^4$,
  and hence $G=\mathcal{G}_2$.

Consider the map $\alpha: \sigma_1 \mapsto \sigma_1^{-1}, \sigma_2 \mapsto \sigma_2^{-1}$. It is easy to check that 
$$(\sigma_1^{-1})^p=(\sigma_{2}^{-1})^{2p^2}=(\sigma_1^{-1}\sigma_2^{-1})^2=[\sigma_1^{-1}, ((\sigma_2^{-1})^2)]^p=[\sigma_1^{-1}, (\sigma_2^{-1})^2, \sigma_1^{-1}]=[\sigma_1^{-1}, (\sigma_2^{-1})^2, (\sigma_2^{-1})^2]=1,$$
meaning $\alpha \in \Aut(G)$, which contradicts the hypothesis that $\mathcal{P}$ is chiral.
  
\noindent{\bf Subcase 1.3:} $\{k_1, k_2\}=\{p^2, 2p\}$.
Then $o(\sigma_1)=p^2$ and $o(\sigma_2^2)=p$. In a similar way, we obtain $G=\mathcal{G}_3$, where
$$\mathcal{G}_3=\lg \sigma_1, \sigma_2  \ | \ \sigma_1^{p^2}, \sigma_2^{2p}, (\sigma_1\sigma_2)^2, [\sigma_1, \sigma_2^2]^p, [\sigma_1, \sigma_2^2, \sigma_1], [\sigma_1, \sigma_2^2, \sigma_2^2] \rg.$$
But this is impossible because $\alpha \in \Aut(G)$, where $\alpha$ taking $(\sigma_1, \sigma_2)$ to $(\sigma_1^{-1}, \sigma_2^{-1})$, contradicting the chirality of $\mathcal{P}$. 

\noindent {\bf Case 2:} $c(P)=3$. Then since ${\rm exp}(P') = p$, $$[\sigma_1, \sigma_2^2]^p=[\sigma_1, \sigma_2^2, \sigma_1]^p=[\sigma_1, \sigma_2^2, \sigma_2^2]^p=1$$
\noindent and as $c(P)=3$,
$$[\sigma_1, \sigma_2^2, \sigma_1, \sigma_1]=[\sigma_1, \sigma_2^2, \sigma_1, \sigma_2^2]=[\sigma_1, \sigma_2^2, \sigma_2^2, \sigma_1]=[\sigma_1, \sigma_2^2, \sigma_2^2, \sigma_2^2]=1.$$

 \noindent{\bf Subcase 2.1:} $\{k_1, k_2\}=\{p, 2p\}$. Then $o(\sigma_1)=o(\sigma_2^2)=p$. Since $P' \cong \mathbb{Z}_p \times \mathbb{Z}_p$, we have 
 $$[\sigma_1, \sigma_2^2, \sigma_1]^i[\sigma_1, \sigma_2^2, \sigma_2^2]^j=1$$
for some $i, j$. This is also enough to define the entire group $G$. We claim that $G=\mathcal{G}_4$, where
 
$$\mathcal{G}_4= \lg \sigma_1, \sigma_2  \ | \ \sigma_1^p, \sigma_2^{2p}, (\sigma_1\sigma_2)^2,
[\sigma_1, \sigma_2^2]^p, [\sigma_1, \sigma_2^2, \sigma_1]^p, [\sigma_1, \sigma_2^2, \sigma_2^2]^p,$$
$$ [\sigma_1, \sigma_2^2, \sigma_1]^i[\sigma_1, \sigma_2^2, \sigma_2^2]^j, [\sigma_1, \sigma_2^2, \sigma_1, \sigma_1], [\sigma_1, \sigma_2^2, \sigma_1, \sigma_2^2]\rg.$$
Obviously, $(G, \{\sigma_1, \sigma_2\})$ satisfies all of the defining relations of $(\mathcal{G}_4, \{\sigma_1, \sigma_2\})$, and hence $G$ is a quotient group of $\mathcal{G}_4$. 
On the other hand, 
consider the quotient group 
$\overline{\mathcal{G}_4}=\mathcal{G}_4/\lg \sigma_1, \sigma_2^2\rg=\lg \overline{\sigma_2} \ | \  (\overline{\sigma_2})^2 \rg$. Then $|\overline{\mathcal{G}_4}| \leq 2$.
For the subgroup $\lg \sigma_1, \sigma_2^2\rg$ of $\mathcal{G}_4$, since $[\sigma_1, \sigma_2^2, \sigma_1, \sigma_1]=[\sigma_1, \sigma_2^2, \sigma_1, \sigma_2^2]=1$, we have $[\sigma_1, \sigma_2^2, \sigma_1] \in Z(\lg \sigma_1, \sigma_2^2\rg)$. 
Let $\widetilde{P}=\lg \sigma_1, \sigma_2^2\rg/\lg [\sigma_1, \sigma_2^2, \sigma_1] \rg=\lg \widetilde{\sigma_1}, \widetilde{\sigma_2^2}\rg$. Then $|\lg \sigma_1, \sigma_2^2\rg| \leq |\widetilde{P}| \cdot p$.
Furthermore, since $\widetilde{P}'=\lg [\widetilde{\sigma_1}, \widetilde{\sigma_2^2}]\rg$, we have
$$|\widetilde{P}|=|\lg \widetilde{\sigma_1}, \widetilde{\sigma_2^2}\rg / \lg [\widetilde{\sigma_1}, \widetilde{\sigma_2^2}]\rg| \cdot |\lg [\widetilde{\sigma_1}, \widetilde{\sigma_2^2}]\rg| \leq p \cdot p \cdot p=p^3.$$
It follows that $|\mathcal{G}_4|\leq 2p^4$, and hence $G=\mathcal{G}_4$.

Consider the map $\alpha: \sigma_1 \mapsto \sigma_1^{-1}, \sigma_2 \mapsto \sigma_2^{-1}$. It is easy to check that 
$$(\sigma_1^{-1})^p=(\sigma_{2}^{-1})^{2p^2}=(\sigma_1^{-1}\sigma_2^{-1})^2=1,$$
$$[\sigma_1^{-1}, ((\sigma_2^{-1})^2)]^p=[\sigma_1^{-1}, (\sigma_2^{-1})^2, \sigma_1^{-1}]^p=[\sigma_1^{-1}, (\sigma_2^{-1})^2, (\sigma_2^{-1})^2]^p=1,$$
$$[\sigma_1^{-1}, (\sigma_2^{-1})^2, \sigma_1^{-1}, \sigma_1^{-1}]=[\sigma_1^{-1}, (\sigma_2^{-1})^2, \sigma_1^{-1}, (\sigma_2^{-1})^2]=1.$$
\noindent Furthermore, by Proposition~\ref{commutator}, 
$$ [\sigma_1^{-1}, (\sigma_2^{-1})^2, \sigma_1^{-1}]^i[\sigma_1^{-1}, (\sigma_2^{-1})^2, (\sigma_2^{-1})^2]^j=1.$$
\noindent Hence $\alpha \in \Aut(G)$, which contradicts the fact that $\mathcal{P}$ is chiral.

 \noindent{\bf Subcase 2.2:} $\{k_1, k_2\}=\{p, 2p^2\}$. Then $o(\sigma_1)=p, o(\sigma_2^2)=p^2$. Moreover, 
 $$\exp(P)=p^2 \mbox{and} \ \mho_1(P)=\lg \sigma_2^{2p}\rg.$$
 Since $P_3= \mho_1(P)$, we have 
 $[\sigma_1, \sigma_2^2, \sigma_1]=(\sigma_2^{2p})^i$ and $[\sigma_1, \sigma_2^2, \sigma_2^2]=(\sigma_2^{2p})^j$ for some $i, j$.
 This is enough to define the entire group $G$. We claim that $G=\mathcal{G}_5$, where
 
$$\mathcal{G}_5= \lg \sigma_1, \sigma_2  \ | \ \sigma_1^p, \sigma_2^{2p^2}, (\sigma_1\sigma_2)^2,
[\sigma_1, \sigma_2^2]^p, [\sigma_1, \sigma_2^2, \sigma_1]^p, [\sigma_1, \sigma_2^2, \sigma_2^2]^p,[\sigma_1, \sigma_2^2, \sigma_1]=(\sigma_2^{2p})^i, $$

$$ 
[\sigma_1, \sigma_2^2, \sigma_2^2]=(\sigma_2^{2p})^j,
 [\sigma_1, \sigma_2^2, \sigma_1, \sigma_1], [\sigma_1, \sigma_2^2, \sigma_1, \sigma_2^2]\rg.$$
Obviously, $(G, \{\sigma_1, \sigma_2\})$ satisfies all of the defining relations of $(\mathcal{G}_5, \{\sigma_1, \sigma_2\})$, and hence $G$ is the quotient group of $\mathcal{G}_5$. Moreover, the same argument as used in Case 2.1 shows that $|\mathcal{G}_5| \leq 2p^4$. Thus, $G=\mathcal{G}_5$.

Consider the map $\alpha: \sigma_1 \mapsto \sigma_1^{-1}, \sigma_2 \mapsto \sigma_2^{-1}$. It is easy to check that 
$$(\sigma_1^{-1})^p=(\sigma_{2}^{-1})^{2p^2}=(\sigma_1^{-1}\sigma_2^{-1})^2=1,$$
$$[\sigma_1^{-1}, ((\sigma_2^{-1})^2)]^p=[\sigma_1^{-1}, (\sigma_2^{-1})^2, \sigma_1^{-1}]^p=[\sigma_1^{-1}, (\sigma_2^{-1})^2, (\sigma_2^{-1})^2]^p=1,$$
$$[\sigma_1^{-1}, (\sigma_2^{-1})^2, \sigma_1^{-1}, \sigma_1^{-1}]=[\sigma_1^{-1}, (\sigma_2^{-1})^2, \sigma_1^{-1}, (\sigma_2^{-1})^2]=1.$$

\noindent Furthermore, by Proposition~\ref{commutator}, 
$$ [\sigma_1^{-1}, (\sigma_2^{-1})^2, \sigma_1^{-1}]=[\sigma_1, \sigma_2^2, \sigma_1]^{-1}=(\sigma_2^{2p i})^{-1}=(\sigma_2^{-1})^{2pi}, $$
$$ [\sigma_1^{-1}, (\sigma_2^{-1})^2, (\sigma_2^{-1})^2]=[\sigma_1, \sigma_2^2, \sigma_2^2]^{-1}=(\sigma_2^{2pj})^{-1}=(\sigma_2^{-1})^{2pj}.$$
\noindent Hence $\alpha \in \Aut(G)$, which contradicts the fact that $\mathcal{P}$ is chiral. 

 \noindent{\bf Subcase 2.3:} $\{k_1, k_2\}=\{p^2, 2p\}$. Then $o(\sigma_1)=p^2, o(\sigma_2^2)=p$. Moreover, 
 $$\exp(P)=p^2 \mbox{and} \ \mho_1(P)=\lg \sigma_1^{p}\rg.$$
 Since $P_3= \mho_1(P)$, we have 
 $[\sigma_1, \sigma_2^2, \sigma_1]=\sigma_1^{pi}$ and $[\sigma_1, \sigma_2^2, \sigma_2^2]=\sigma_1^{pi}$ for some $i, j$.
 In a similar way, we obtain $G=\mathcal{G}_6$, where
$$\mathcal{G}_6= \lg \sigma_1, \sigma_2  \ | \ \sigma_1^{p^2}, \sigma_2^{2p}, (\sigma_1\sigma_2)^2,
[\sigma_1, \sigma_2^2]^p, [\sigma_1, \sigma_2^2, \sigma_1]^p, [\sigma_1, \sigma_2^2, \sigma_2^2]^p, [\sigma_1, \sigma_2^2, \sigma_1]=\sigma_1^{pi}, $$
$$ 
[\sigma_1, \sigma_2^2, \sigma_2^2]=\sigma_1^{pj},
 [\sigma_1, \sigma_2^2, \sigma_1, \sigma_1], [\sigma_1, \sigma_2^2, \sigma_1, \sigma_2^2]\rg.$$
But this will again give an element $\alpha \in \Aut(G)$, with $\alpha(\sigma_1, \sigma_2)=(\sigma_1^{-1}, \sigma_2^{-1})$, contradicting the chirality of $\mathcal{P}$. \qed

\section{Proof of Theorem~\ref{maintheorem3}}\label{sec6}

The main purpose of this section is to prove Theorem~\ref{maintheorem3}.
In order to do so, we will give a constructive proof.  Let $P, \sigma, \tau$ be the group and maps defined in Section~\ref{sec3}. We have two cases to consider, namely the case where $m$ is odd and at least 5 and the case where $m$ is even and at least $p+3$.

\noindent {\bf Case 1:} $m$ is odd and $m \geq 5$.

Let $G=P \rtimes \lg \sigma \rg$. Then $|G|=2p^m$ and $G=\lg s_1, \b, \sigma \rg$. Let $\sigma_1=\beta s_1, \sigma_2=s_1^{-1}\sigma$. We 
will show that $(G, \{\sigma_1, \sigma_2\})$  
is the automorphism group of a chiral polytope of type $\{p^2, 2p^e\}$.

\noindent{\bf Claim 1:} $G=\lg \sigma_1, \sigma_2\rg$, and $[o(\sigma_1), o(\sigma_2), o(\sigma_1\sigma_2)]=[p^2, 2p^e, 2]$.

Obviously, $\sigma_2^2=(s_1^{-1}\sigma)^2=s_1^{-2}$. Since $o(s_1)=p^e$, 
$s_1 \in \lg s_1^{-2}\rg =\lg \sigma_2^2\rg \leq \lg \s_1, \s_2\rg$. It follows that $\s_1 s_1^{-1}=\beta \in \lg \s_1, \s_2\rg$, 
$s_1 \s_2=\sigma \in \lg \s_1, \s_2\rg$ and hence $\lg \s_1, \s_2\rg=\lg s_1, \b, \sigma \rg=G$.

Next, since $o(\sigma)=2$ and $\b^{\sigma}=\b^{-1}$, we have
$$(\sigma_1\sigma_2)^2=(\beta s_1s_1^{-1}\sigma)^2=(\beta \sigma)^2=\beta\beta^{\sigma}=\beta\beta^{-1}=1,$$
and hence $o(\sigma_1\sigma_2)=2$.
Moreover, since $\sigma_1 \in P  \backslash  A$ and $\sigma_2^2=s_1^{-2}$, we have $o(\sigma_1)=p^2$ and $o(\sigma_2)=2p^e$, as claimed.

\noindent{\bf Claim 2:} $\lg \sigma_1\rg \cap \lg \sigma_2\rg=\{1\}$.

Since $o(\sigma_1)=p^2$ and $\lg \sigma_1 \rg \cap \lg \sigma_2\rg  \leq \lg \sigma_1 \rg$, $|\lg \sigma_1 \rg \cap \lg \sigma_2\rg|=1, p$ or $p^2$.
 If $|\lg \sigma_1 \rg \cap \lg \sigma_2\rg|=p^2$, then $\lg \sigma_1 \rg \cap \lg \sigma_2\rg=\lg \sigma_1\rg$. It follows that $\sigma_1 \in \lg \sigma_2\rg$ 
and hence $G=\lg \sigma_1, \sigma_2\rg =\lg \sigma_2 \rg$, which is impossible because $G$ is non-abelian.

If $|\lg \sigma_1 \rg \cap \lg \sigma_2\rg|=p$, then $\lg \sigma_1\rg \cap \lg \sigma_2 \rg=\lg \sigma_1^p\rg =\lg \sigma_2^{2p^{e-1}}\rg$.
Observe that $\sigma_1^p=(\beta s_1)^p=\beta^p=s_r^{p^{e-1}}$ and $\lg \sigma_2^{2p^{e-1}}\rg=\lg (s_1^{-2})^{p^{e-1}}\rg=\lg s_1^{p^{e-1}}\rg$.
However, by Lemma~\ref{reven}, $r$ is even, thus $s_r^{p^{e-1}} \notin \lg s_1^{p^{e-1}}\rg$. It means that $\lg \sigma_1^p\rg \ne \lg \sigma_2^{2p^{e-1}}\rg$, and 
hence $|\lg \sigma_1 \rg \cap \lg \sigma_2\rg| \ne p$. Thus, $\lg \sigma_1\rg \cap \lg \sigma_2\rg=\{1\}$, as claimed.

\noindent{\bf Claim 3:} There is no $\alpha \in \Aut(G)$ such that  $\alpha(\sigma_1, \sigma_2) = (\sigma_1^{-1}, \sigma_2^{-1})$. 

Suppose that there exists such an $\alpha$. 
Since $P$ is a characteristic subgroup of $G$, the automorphism $\alpha$ of $G$ induces an automorphism of $P$, say $\alpha|_{P}$. 
Note that  $\sigma_2^2=s_1^{-2}$, so $(s_1^{-2})^{\alpha}=((\sigma_2)^2)^{\alpha}=((\sigma_2)^{\alpha})^2=\sigma_2^{-2}=s_1^2$. It follows that $s_1^{\alpha}=s_1^{-1}$.
Moreover, since $s_1^{-1}\beta^{-1}=(\beta s_1)^{-1}=\sigma_1^{-1}=(\sigma_1)^{\alpha}=(\beta\s_1)^{\alpha}=\beta^{\alpha}s_1^{\alpha}=\beta^{\alpha}s_1^{-1}$, we have
$\beta^{\alpha}=s_1^{-1}\beta^{-1}s_1$. Thus,  $\alpha|_{P}$ takes $(s_1, \beta)$ to $({s_1^{-1}, (\beta^{s_1})^{-1}})$. 
Let $\sigma(s_1)$ be the inner automorphism of $P$ induced by $s_1$.
Then we have the following sequence: 
$$\alpha|_{P} \in \Aut(P) \longrightarrow \alpha|_{P} \circ \sigma(s_1) \in \Aut(P) \longrightarrow \alpha|_{P} \circ \sigma(s_1)\circ \sigma \in \Aut(P) \Leftrightarrow \tau \in \Aut(P)$$

\noindent However, by Lemma~3.1, $\tau \notin \Aut(P)$, a contradiction, and hence $\alpha \notin \Aut(G)$, as claimed.

Claims 1, 2 and 3 show that $(G, \{\sigma_1, \sigma_2\})$  
is the automorphism group of a chiral polytope of type $\{p^2, 2p^e\}$. 

\noindent{\bf Case 2:} $m$ is even, and $m \geq p+3$. 

We define the group $$G^{\ast}=\lg s_1^{\ast}, s_2^{\ast}, \cdots, s_r^{\ast}, \cdots, s_{p-1}^{\ast}, \beta^{\ast}, \sigma^{\ast}, x \rg$$ with the following relations:

$$ s_1^{\ast}{}^{p^e}, s_2^{\ast}{}^{p^e}, \cdots, s_r^{\ast}{}^{p^e}, s_{r+1}^{\ast}{}^{p^{e-1}}, \cdots s_{p-1}^{\ast}{}^{p^{e-1}}, $$
$$\b^{\ast}{}^{p^2}, \b^{\ast}{}^{p}=s_r^{\ast}{}^{p^{e-1}}, x^p,$$
$$[s_i^{\ast}, s_j^{\ast}],  \ {\rm for} \ 1 \leq i, j \leq p-1,\ \mbox{and}\ s_{k+1}^{\ast}=[s_k^{\ast}, \b^{\ast}],  \ {\rm for} \ 1 \leq k \leq p-2,$$
$$x=s_{1}^{\ast} {}^{-\tbinom{p}{1}}s_2^{\ast} {}^{-\tbinom{p}{2}}\cdots s_{p-2}^{\ast}{}^{-\tbinom{p}{p-2}}s_{p-1}^{\ast}{}^{-\tbinom{p}{p-1}}[s_{p-1}^{\ast}, \b^{\ast}]^{-1},$$
$$s_1^{\ast}{}^{\sigma^{\ast}}=s_1^{\ast}, \beta^{\ast}{}^{\sigma^{\ast}}=\beta^{\ast}{}^{-1}, [s_1^{\ast}, x]=[\beta^{\ast}, x]=[\sigma^{\ast},x]=1.$$
where $1 \leq e,1 \leq r \leq p-1$, and $r$ is even. 

Obviously, $G^{\ast}=\lg s_1^{\ast}, \b^{\ast}, \sigma^{\ast}\rg$. Let  $\sigma_1^{\ast}=\beta^{\ast} s_1^{\ast},  \sigma_2^{\ast}=s_1^{\ast}{}^{-1}\sigma^{\ast}$, and $e \geq 2$.

\noindent{\bf Claim 4:} $o(x)=p$ and $|G^{\ast}|=2p^{(p-1)(e-1)+r+2} \geq 2p^{p+3}$.

Since $[s_1^{\ast}, x]=[\beta^{\ast}, x]=[\sigma^{\ast}, x]=1$, we have $x \in Z(G^{\ast})$. Consider the quotient group $G^{\ast}/\lg x\rg$.
Then $G^{\ast}/\lg x\rg \cong G$ (which is the group defined in Case 1) and hence $|G^{\ast}/\lg x\rg|=2p^{(p-1)(e-1)+r+1}$. We only need to prove that $o(x) =p$. 
Then $G^{\ast}$ is a central extension of $G$ by $\lg x \rg$.
Unfortunately, $x$ not always has order $p$, which is why we need $e \geq 2$.

Let $H=\lg a, b \ | \ a^{p^2}, b^p, [a, b] \rg \cong \mathbb{Z}_{p^2} \times \mathbb{Z}_{p}$.
It is easy to check that $(H, \{a, b, 1\})$ satisfies all of the relations of $(G^{\ast}, \{s_1^{\ast}, \beta^{\ast}, \sigma^{\ast}\})$.
Hence the map: $(s_1^{\ast}, \beta^{\ast}, \sigma^{\ast}) \mapsto (a, b, 1)$ induces an epimorphism $\pi$ from $G^{\ast}$ to $H$.
Moreover, $\pi :x \mapsto a^{-p}$. Note that $o(a) =p^2$, so $o(\pi(x))=p$. It follows that $o(x)=p$ and hence 
$|G^{\ast}|=p|G|=2p^{(p-1)(e-1)+r+2}$.

From now on, let $|G^{\ast}|=2p^m$. Since $r$ is even and $e \geq 2$, we get $|G^{\ast}|=2p^{(p-1)(e-1)+r+2} \geq 2p^{p+3}$ which is fine as $m$ is even 
and $m \geq p+3$ by hypothesis.

\noindent{\bf Claim 5:}  $G^{\ast}=\lg \sigma_1^{\ast}, \sigma_2^{\ast}\rg$, and  $[o(\sigma_1^{\ast}), o(\sigma_2^{\ast}), o(\sigma_1^{\ast}\sigma_2^{\ast})]=[p^2, 2p^e, 2]$.

The same argument as used in Case 1 shows that $G^{\ast}=\lg \sigma_1^{\ast}, \sigma_2^{\ast}\rg$ and $o(\sigma_2^{\ast})=2p^e, o(\sigma_1^{\ast}\sigma_2^{\ast})=2$.
Consider the quotient group $G^{\ast} / \lg x\rg \cong G$. 
Note that $\sigma_1^p = \b^p$ in $G$ implies that $\sigma_1^{\ast p} =\b^{\ast p} x^i$ in $G^{\ast}$ for some $i$. 
It follows that $\sigma_1^{\ast p^2}=(\b^{\ast p}x^i)^p=\b^{\ast p^2}(x^i)^p=1$ in $G^{\ast}$. Then $o(\sigma_1^{\ast})=p^2$.

\noindent{\bf Claim 6:} $\lg \sigma_1^{\ast}\rg \cap \lg \sigma_2^{\ast}\rg=\{1\}$.

Let $\gamma$ be the epimorphism from $(G^{\ast},\{\sigma_1^{\ast}, \sigma_2^{\ast}\})$ onto $(G, \{\sigma_1, \sigma_2\})$. 
Note that $o(\sigma_1^{\ast})=o(\sigma_1)=p^2$ and $\lg \sigma_1 \rg \cap \lg \sigma_2\rg=\{1\}$. Then Proposition~\ref{quotient criterion} implies that $\lg \sigma_1^{\ast}\rg \cap \lg \sigma_2^{\ast}\rg=\{1\}$.

\noindent{\bf Claim 7:} There is no $\alpha^\ast \in \Aut(G^\ast)$ such that $\alpha^{\ast}(\sigma_1^{\ast}, \sigma_2^{\ast})=(\sigma_1^{\ast -1}, \sigma_2^{\ast -1})$. 

Suppose that such an $\alpha^{\ast}$ exists.
Then the same argument as used in Case 1 shows that $\alpha^{\ast}|_{P^{\ast}}$ takes $(s_1^{\ast}, \beta^{\ast})$ to $({s_1^{\ast -1}, (\beta^{\ast s_1^{\ast}})^{-1}})$. 
On the other hand, 
since $| \lg \sigma^{\ast}\rg|=2$ and $|P^{\ast}|=p^m$, we get $P^{\ast} \cap \lg \sigma^{\ast}\rg=\{1\}$.
It follows that $G^{\ast}=P^{\ast} \rtimes \lg \sigma^{\ast}\rg$, and hence $\sigma^{\ast} \in \Aut(P^{\ast})$.
Let $\tau^{\ast}= \alpha^{\ast}|_{P^{\ast}} \circ \sigma(s_1^{\ast})\circ \sigma^{\ast}$. Then $(s_1^{\ast}, \b^{\ast})^{\tau^{\ast}}=(s_1^{\ast -1}, \b^{\ast})$
 and $\tau^{\ast} \in \Aut(P^{\ast})$.
However, the same argument as used in Lemma 3.1 shows that $\tau^{\ast} \notin \Aut(P^{\ast})$, a contradiction.

Claims 4 to 7 above imply that $(G^{\ast}, \{\sigma_1^{\ast}, \sigma_2^{\ast}\})$  is the automorphism group of a chiral polytope of type $\{p^2, 2p^e\}$.  \qed

\section{Acknowledgements} This work was supported by the National Natural Science Foundation of China (12201371,
12271318,12331013,12311530692,12271024,12161141005), the 111 Project of China (B16002), an Action de Recherche Concertée grant of the Communauté Française Wallonie Bruxelles and a PINT-BILAT-M grant from the Fonds National de la Recherche Scientifique de Belgique (FRS-FNRS).

\end{document}